\documentclass[11pt]{tran-l}

\newtheorem{theorem}{Theorem}[section]
\newtheorem{lemma}[theorem]{Lemma}
\newtheorem{prop}[theorem]{Proposition}
\newtheorem{conj}[theorem]{Conjecture}
\newtheorem{cor}[theorem]{Corollary}

\theoremstyle{definition}

\theoremstyle{remark}
\newtheorem{remark}[theorem]{Remark}

\numberwithin{equation}{section}

\begin{document}

\title[Relative Completions]{Relative Completions of Linear Groups over ${\mathbb Z}[t]$ and ${\mathbb Z}[t,t^{-1}]$}

% author one information
\author{Kevin P. Knudson}
\address{Department of Mathematics, Northwestern University, Evanston, IL 60208}
\email{knudson@math.nwu.edu}
\thanks{Supported by an NSF Postdoctoral Fellowship, grant no. DMS-9627503}

\newcommand{\powz}{SL_n({\mathbb Z}[[T]])}
\newcommand{\slnz}{SL_n({\mathbb Z})}
\newcommand{\zt}{SL_n({\mathbb Z}[t])}
\newcommand{\zl}{SL_n({\mathbb Z}[t,t^{-1}])}
\newcommand{\slnq}{SL_n({\mathbb Q})}
\newcommand{\powq}{SL_n({\mathbb Q}[[T]])}
\newcommand{\zz}{{\mathbb Z}}
\newcommand{\zq}{{\mathbb Q}}
\newcommand{\U}{{\mathcal U}}
\newcommand{\cP}{{\mathcal P}}
\newcommand{\M}{{\mathcal M}}
\newcommand{\G}{{\mathcal G}}
\newcommand{\E}{{\mathcal E}}
\newcommand{\lra}{\longrightarrow}
\newcommand{\ra}{\rightarrow}
\newcommand{\Lra}{\Longrightarrow}
\newcommand{\slnr}{SL_n(R)}
\newcommand{\lslnz}{{\mathfrak sl}_n({\mathbb Z})}
\newcommand{\lslnq}{{\mathfrak sl}_n({\mathbb Q})}
\newcommand{\bop}{\bigoplus}

\subjclass{Primary 55P60, 20G35, 20H05; Secondary 20G10, 20F14}
\date{\today}

\begin{abstract}We compute the completion of the groups $SL_n({\mathbb Z}[t])$ and $SL_n({\mathbb Z}[t,t^{-1}])$ relative to the obvious homomorphisms to $SL_n({\mathbb Q})$; this is a generalization of the classical      Malcev completion.  We also make partial computations of the rational second cohomology of these groups. 
\end{abstract}

\maketitle

The Malcev (or $\zq$-) completion of a group $\Gamma$ is a prounipotent group $\cP$ defined over $\zq$ together
with a homomorphism $\varphi:\Gamma \ra \cP$ satisfying the following universal mapping property:  If
$\psi:\Gamma \ra \U$ is a map of $\Gamma$ into a prounipotent group then there is a unique map
$\Phi:\cP \ra \U$ such that $\psi = \Phi\varphi$.  If $H_1(\Gamma,\zq) = 0$, then the group $\cP$ is trivial
and is therefore useless for studying $\Gamma$.  In particular, the Malcev completions of the groups
$\zt$ and $\zl$ are trivial when $n\ge 3$ (this follows from the work of Suslin \cite{suslin}).

Here we consider Deligne's notion of relative completion.  Suppose $\rho:\Gamma \ra S$ is a
representation of $\Gamma$ in a semisimple linear algebraic group over $\zq$.  Suppose that the image of
$\rho$ is Zariski dense in $S$.  The completion of $\Gamma$ relative to $\rho$ is a proalgebraic group
$\G$ over $\zq$, which is an extension of $S$ by a prounipotent group $\U$, and homomorphism
$\tilde{\rho}:\Gamma \ra \G$ which lifts $\rho$ and has Zariski dense image.  When $S$ is the trivial group,
$\G$ is simply the classical Malcev completion.  The relative completion satisfies an obvious universal
mapping property.  The basic theory of relative completion was developed by R. Hain \cite{hain2} (and
independently by E. Looijenga (unpublished)), and is reviewed in Section \ref{relcomp} below.

In this paper we consider the completions of the groups $\zt$ and $\zl$ relative to the homomorphisms to
$\slnq$ given by setting $t=0$ (respectively, $t=1$).  There is an obvious candidate for the relative
completion, namely the proalgebraic group $\powq$.  The map $$\zt \ra \powq$$ is the obvious inclusion and the
map $$\zl \ra \powq$$ is induced by the map $t\mapsto 1+T$.

\medskip

\noindent {\bf Theorem.}  {\em  For all $n\ge 3$, the group $\powq$ is the relative completion of both
$\zt$ and $\zl$.}

\medskip

\noindent {\em Remark.}  We expect that the theorem holds for an arbitrary simple group $G$ of sufficiently
large rank (large enough to guarantee the vanishing of $H^2(G(\zz),A)$ for nontrivial $G(\zq)$-modules $A$).
We have chosen to work with $SL_n$ just to be concrete.

\medskip

The theorem does not hold for $n=2$ (see Section \ref{n2} below).  Our proof breaks down in this case essentially because the
$\zz$-Lie algebra ${\mathfrak sl}_2(\zz)$ is not perfect.

We use this result to study the cohomology of the groups $\zt$ and $\zl$.  
This is motivated by an attempt to find unstable analogues of the Fundamental
Theorem of Algebraic $K$--theory. Recall that if $A$ is a regular ring, then
there are natural isomorphisms $K_\bullet(A[t]) \cong K_\bullet(A)$ and
$K_\bullet(A[t,t^{-1}]) \cong K_\bullet(A) \oplus K_{\bullet - 1}(A)$.
An unstable analogue does exist for infinite fields:  If $k$ is an 
infinite field, then $H_\bullet(SL_n(k[t]),\zz) \cong H_\bullet(SL_n(k),\zz)$
for all $n$ \cite{knudson}.  Since $\zz$ is regular, one might guess that
an analogous statement holds for $n$ sufficiently large.  We note, however,
that if such a result holds, we must have $n\ge 3$ since $H_1(SL_2(\zz[t]),\zz)$
has infinite rank \cite{gmv}, while $H_1(SL_2(\zz),\zz) \cong \zz/12$.

  The basic idea is to use continuous 
cohomology.  Following Hain \cite{hain1}, we define the continuous cohomology of a group $\pi$ to be
$$H^\bullet_{\text{cts}}(\pi,\zq) = \varinjlim H^\bullet(\pi/\Gamma^r\pi,\zq),$$
where $\Gamma^\bullet \pi$ denotes the lower central series of $\pi$.  There is a natural map
$$H^\bullet_{\text{cts}}(\pi,\zq) \lra H^\bullet(\pi,\zq)$$
which is injective in degree $2$ provided $H_1(\pi,\zq)$ is finite dimensional.

Consider the extension
$$1 \lra K(R) \lra \slnr \lra \slnz \lra 1$$
for $R=\zz[t], \zz[t,t^{-1}]$.  This yields a spectral sequence for
computing  the rational cohomology of $SL_n(R)$.
In light of the following result it is reasonable to conjecture that
$H^2(\slnr,\zq)=0$ for $n\ge 3$.

\medskip

\noindent {\bf Theorem.} {\em  If $n\ge 3$, then $H^0(\slnz,H^2_{\text{\em cts}}(K(R),\zq))=0$.}

\medskip

Of course, one can see that $H^2(\slnr,\zq)=0$ for $n\ge 5$ by using van der Kallen's stability theorem
\cite{vdk} and the Fundamental Theorem of Algebraic K-theory.  The above result provides evidence for the
vanishing of $H^2(\slnr,\zq)$ for $n=3,4$.  We note, however, that $H^2(SL_2(\zz[t]),\zq)$ has infinite rank
(this is a consequence of results of Grunewald, et. al. \cite{gmv}).

The study of the relative completion of the fundamental group of a
complex algebraic variety $X$ is related to the study of variations of mixed Hodge structure over $X$ 
\cite{hain2}.  Moreover,
relative completions were used with great success by R. Hain in his study of mapping class groups
$\M_g$ and Torelli groups ${\mathcal T}_g$ \cite{hain2},\cite{hain3}.  In particular, he was able to provide
a presentation of the Malcev Lie algebra of ${\mathcal T}_g$ which in turn gives a partial computation of
$H^2({\mathcal T}_g,\zq)$. This also yields a description of the completion
${\mathcal G}_g$ of ${\mathcal M}_g$ with respect to its representation on
the first homology of the surface.  However, the map ${\mathcal M}_g \ra
{\mathcal G}_g$ remains a mystery.    As far as we know, the results of this paper provide the first concrete descriptions of
relative completions and the associated homomorphisms aside from the obvious trivial ones ({\em e.g.}, $\slnz \ra \slnq$, $n\ge 3$).

This paper obviously owes a great deal to the work of Dick Hain and I thank him for suggesting this
problem to me.

\section{Malcev Completions}\label{malcomp}

Recall that the Malcev completion of a group $\Gamma$ is a prounipotent group $\M$ over $\zq$, together with a
map $\Gamma \ra \M$ which satisfies the obvious universal mapping property.  We recall the construction of 
$\M$ as given by Bousfield \cite{bousfield}.

First suppose that $G$ is a nilpotent group.  The Malcev completion of $G$ consists of a group $\widehat{G}$ and
a homomorphism $j:G\ra \widehat{G}$.  It is characterized by the following three properties 
\cite[Appendix A, Cor. 3.8]{quillen2}:
\begin{enumerate}
\item $\widehat{G}$ is nilpotent and uniquely divisible.
\item The kernel of $j$ is the torsion subgroup of $G$.
\item If $x \in \widehat{G}$, then $x^n \in \text{im}\; j$ for some $n\ne 0$.
\end{enumerate}
Quillen constructs $\widehat{G}$ as the set of grouplike elements of the completed group algebra $\widehat{\zq}G$ (completed with respect to the augmentation ideal).

Now, if $G$ is an arbitrary group, denote by $G_r$ the nilpotent group $G/\Gamma^rG$, where $\Gamma^\bullet G$ is
the lower central series of $G$.  Following Bousfield \cite{bousfield}, we define the Malcev completion of $G$ to be
$$\M = \varprojlim \widehat{G}_r$$
where $\widehat{G}_r$ is the Malcev completion of $G_r$.  One checks easily that the group $\M$ satisfies the
universal mapping property.

\section{Relative Completions}\label{relcomp}

In this section we review the theory of relative completion.  All results in this section are due to
R. Hain \cite{hain2}.  

Let $\Gamma$ be a group and $\rho:\Gamma \ra S$ a Zariski dense representation of $\Gamma$ in a semisimple
algebraic group $S$ over $\zq$.  The completion of $\Gamma$ relative to $\rho$ may be constructed as
follows.  Consider all commutative diagrams of the form
$$\begin{array}{ccccccccc}
1 & \ra & U & \ra &     E                 &       \ra      & S & \ra & 1 \\
  &     &   &     & \tilde{\rho} \uparrow & \nearrow \rho  &   &     &   \\
  &     &   &     &  \Gamma               &                &   &     &
\end{array}$$
where $E$ is a linear algebraic group over $\zq$, $U$ is a unipotent subgroup of $E$, and the image of 
$\tilde{\rho}$ is Zariski dense.  The collection of all such diagrams forms an inverse system \cite[Prop 2.1]{hain2}
and we define the completion of $\Gamma$ relative to $\rho$ to be
$$\G = \varprojlim E.$$
The group $\G$ satisfies the following universal mapping property.  Suppose that $\E$ is a proalgebraic group over
$\zq$ such that there is a map $\E \ra S$ with prounipotent kernel.  If $\varphi:\Gamma \ra \E$ is a homomorphism
whose composition with $\E \ra S$ is $\rho$, then there is a unique map $\tau:\G \ra \E$ such that the diagram
$$\begin{array}{ccccc}
       &               &    \G      &          &         \\
       & \nearrow      &            & \searrow &         \\
\Gamma &               & \tau\downarrow &      &   S     \\
       &\varphi\searrow &           & \nearrow &         \\
       &               &     \E     &          & 
\end{array}$$
commutes.

Denote by $L$ the image of $\rho:\Gamma \ra S$ and by $T$ the kernel.  Let $\Gamma \ra \G$ be the completion of
$\Gamma$ relative to $\rho$ and let $\U$ be the prounipotent radical of $\G$.  Consider the commutative diagram
$$\begin{array}{ccccccccc}
1 & \ra & T & \ra &\Gamma & \ra & L & \ra & 1 \\
  &     &\downarrow &  &  \downarrow &  &  \downarrow &  &  \\
1 & \ra & \U & \ra & \G & \ra & S & \ra & 1.
\end{array}$$    
Denote by ${\mathcal T}$ the classical Malcev completion of $T$.  The universal mapping property of ${\mathcal T}$ gives a map $\Phi:{\mathcal T}
\ra \U$ whose composition with $T\ra {\mathcal T}$ is the map $T\ra \U$.  Denote the kernel of $\Phi$ by 
${\mathcal K}$.  We have the following two results.

\begin{prop}[\cite{hain2}, Prop. 4.6] \label{surj}  Suppose that the natural action of $L$ on $H_1(T,\zq)$ 
extends to a rational representation of $S$.  If $H^1(L,A)=0$ for all rational representations $A$ of $S$, then
$\Phi$ is surjective. \hfill $\qed$
\end{prop}

\begin{prop}[\cite{hain2}, Prop. 4.13] \label{inj} Suppose $H_1(T,\zq)$ is finite dimensional and that $H^1(L,A)$
vanishes for all rational representations $A$ of $S$.  Suppose further that $H^2(L,A) = 0$ for all nontrivial
rational representations of $S$.  Then there is a surjective map $H_2(L,\zq) \ra {\mathcal K}$. \hfill $\qed$
\end{prop}

Observe that $H^1(\slnz,A)= 0$ for $n\ge 3$ by Raghunathan's theorem \cite{rag}.  Moreover the second condition
that $H^2(\slnz,A)=0$ for all nontrivial $A$ holds for $n\ge 3$ as well \cite{borel}.\footnote{The result in 
\cite{borel} only implies vanishing for $n\ge 8$.  However, this is easily strengthened to $n\ge 3$; compare with
\cite[Theorem 3.2]{hain3}.}

\section{The Malcev Completion of $K(R)$}\label{kercomp}

Consider the short exact sequences 
$$1 \lra K(R) \lra \slnr \lra \slnz \lra 1$$
for $R = \zz[t], \zz[t,t^{-1}]$, and $n\ge 3$.  In this section we compute the Malcev completion of $K(R)$.

Denote by ${\mathfrak m}_{\zz[t]}$ (resp. ${\mathfrak m}_{\zz[t,t^{-1}]}$) the ideal $(t)$ (resp. $(t-1)$) of
$\zz[t]$ (resp. $\zz[t,t^{-1}]$).  For each $l\ge 1$, define a subgroup $K^l(R)$ by
$$K^l(R) = \{ X \in K(R): X \equiv I \mod {\mathfrak m}_R^l\}.$$
One checks easily that $K^\bullet(R)$ is a descending central series; that is, $$[K^i,K^j] \subseteq K^{i+j}.$$
It follows that for each $i$, $\Gamma^iK \subseteq K^i$.

For each $i\ge 1$, define homomorphisms $\rho_i$, $\sigma_i$ as follows.  If $X\in K^i(\zz[t])$, write
$$X= I + t^iX_i +\cdots + t^mX_m$$
where each $X_j$ is a matrix with integer entries.  Define $$\rho_i:K^i(\zz[t])\ra \lslnz$$ by $\rho_i(X) = X_i$.
Similarly, if $Y \in K^i(\zz[t,t^{-1}])$ we may write
$$Y=I + (t-1)^iY_i \mod (t-1)^{i+1}$$ since $(t^{-1} - 1)^i \equiv (-1)^i (t-1)^i \mod (t-1)^{i+1}$.  Now
define $$\sigma_i:K^i(\zz[t,t^{-1}]) \ra \lslnz$$ by $\sigma_i(Y) = Y_i$.  These maps are well-defined since
the condition $\det Z = 1$ in $K^i(R)$ forces $\text{trace}\; Z_i = 0$.  Moreover, it is easy to see that the
maps $\rho_i$, $\sigma_i$ are surjective group homomorphisms with kernel $K^{i+1}$.  Thus for each $i\ge 1$, we
have $$K^i(R)/K^{i+1}(R) \cong \lslnz.$$

Consider the associated graded $\zz$-Lie algebra
$$\text{Gr}^\bullet K(R) = \bop_{i\ge 1} K^i(R)/K^{i+1}(R).$$
If $n\ge 3$, the Lie algebra $\lslnz$ satisfies $\lslnz = [\lslnz,\lslnz]$.  It follows that the graded algebra
$\text{Gr}^\bullet K(R)$ is generated by $\text{Gr}^1 K(R)$.  The following lemma is easily proved (compare with
\cite[Appendix A, Prop. 3.5]{quillen2}).

\begin{lemma}\label{gen}  Let $G$ be a group with filtration $G=G^1 \supseteq G^2 \supseteq \cdots$.  Then the
associated graded Lie algebra $\text{\em Gr}^\bullet G$ is generated by $\text{\em Gr}^1 G$ if and only if 
$G^r = G^{r+1}\Gamma^r$ for each $r\ge 1$. \hfill $\qed$
\end{lemma}

\begin{cor}\label{compcor}  Suppose $\bigcap G^r = \{1\}$.  If $\text{\em Gr}^\bullet G$ is generated by          $\text{\em Gr}^1 G$,
then the completions of $G$ with respect to the filtration $G^\bullet$ and the lower central series
$\Gamma^\bullet G$ are isomorphic; that is,
$$\varprojlim G/G^r \cong \varprojlim G/\Gamma^rG.$$
\end{cor}

\begin{proof}  Consider the short exact sequence
$$1 \lra G^r/\Gamma^r \lra G/\Gamma^r \lra G/G^r \lra 1.$$
Since $\text{Gr}^\bullet G$ is generated by $\text{Gr}^1 G$, we have $G^r = G^{r+1}\Gamma^r$ for each $r$.  It 
follows that the inverse system $\{G^r/\Gamma^r\}$ is surjective.  This, in turn, implies that the natural map
$$\varprojlim G/\Gamma^r \lra \varprojlim G/G^r$$
is surjective.  Injectivity follows since the assumption that $\bigcap G^r = \{1\}$ implies that
$\varprojlim G^r/\Gamma^r = \{1\}$.
\end{proof}

We now compute the Malcev completions of the groups $K(R)/K^i(R)$.  We first provide the following
result.

\begin{lemma}\label{invlim}  The completion of $\zz[t,t^{-1}]$ with respect to the ideal $(t-1)$ is the power
series ring $\zz[[T]]$.  The canonical map $\zz[t,t^{-1}] \ra \zz[[T]]$ sends $t$ to $1+T$.
\end{lemma}

\begin{proof}  This follows easily once one notes that in $\zz[t,t^{-1}]/(t-1)^m$, we have
$t^{-1} = 1 + (t-1) + \cdots + (t-1)^{m-1}$, so that any polynomial in $\zz[t,t^{-1}]/(t-1)^m$ may be written
as a polynomial in nonnegative powers of $(t-1)$.
\end{proof}

Consider the short exact sequence
$$1 \lra \overline{K} \lra \powz \stackrel{T=0}{\lra} \slnz \lra 1.$$

\begin{cor}\label{cor1}  The group $\overline{K}$ is the completion of $K(R)$ with respect to the filtration
$K^1(R) \supset K^2(R) \supset \cdots$ and with respect to the lower central series of $K(R)$.
\end{cor}

\begin{proof}  The first assertion follows from Lemma \ref{invlim} and the second from Corollary \ref{compcor}.
\end{proof}

Observe that the group $\overline{K}$ has a filtration given by powers of $T$ (exactly as $K(\zz[t])$ does) and
that the successive graded quotients are isomorphic to $\lslnz$.  Denote the filtration by $\overline{K}^\bullet$.

We have an analogous sequence over $\zq$:
$$1 \lra \U \lra \powq \stackrel{T=0}{\lra} \slnq \lra 1,$$
and the corresponding $T$-adic filtration $\U^\bullet$ in $\U$.  In this case, the successive graded quotients are
isomorphic to $\lslnq$.  We can assemble our exact sequences into a commutative diagram
$$\begin{array}{ccccccccc}
1  &  \ra & K(R)  & \ra &  \slnr &  \ra & \slnz & \ra & 1 \\
   &       &\downarrow &  & \downarrow &  &  || &  &           \\
1  & \ra  &\overline{K} & \ra & \powz & \stackrel{T=0}{\ra} & \slnz & \ra & 1 \\
   &       &\downarrow &  & \downarrow &  & \downarrow &  &    \\
1  & \ra  & \U  &  \ra & \powq & \stackrel{T=0}{\ra} & \slnq & \ra & 1.
\end{array}$$

\begin{prop}  The map $K(R)/K^r(R) \stackrel{j}{\ra} \U/\U^r$ is the Malcev completion.
\end{prop}

\begin{proof}  According to Quillen's criterion (see Section \ref{malcomp}) we must check three things.  First,
the group $\U/\U^r$ is nilpotent and uniquely divisible.  Nilpotency is obvious, so suppose
$$Y = I + TY_1 + \cdots + T^{r-1}Y_{r-1}$$ is an element of $\U/\U^r$.  For each $n>0$, we must find a unique
$X \in \U/\U^r$ with $X^n = Y$.  For an arbitrary $X \in \U/\U^r$, write $X = I + TX_1 + \cdots + T^{r-1}X_{r-1}$,
and consider the equation
\begin{eqnarray*}
Y &  = &  X^n \\
  &= & I +TnX_1 + T^2(nX_2 + \binom {n}{2}X_1^2) + \cdots + T^{r-1}(nX_{r-1} + p(\{X_i\}_{i=1}^{r-2}))
\end{eqnarray*}
where $p(X_1,\dots ,X_{r-2})$ is a polynomial in the $X_i$, $i \le r-2$.  Clearly, we can solve this equation
inductively for the $X_i$ and find a unique $X$.

Second, we must show that the kernel of $j$ is the torsion subgroup of $K(R)/K^r(R)$.  This is clear since 
$K(R)/K^r(R)$ is torsion-free ({\em i.e.,} if some power of $X\in K$ lies in $K^r$, then $X\in K^r$ already)
and the map is injective.

Finally, we must show that if $X\in \U/\U^r$, then $X^m \in \text{im}\; j$ for some $m\ne 0$.  We prove this by
induction on $r$, beginning at $r=2$.  Let $X=I + TX_1$ be an element of $\U/\U^2$.  Then there is an $m>0$ such
that $mX_1$ consists of integer entries.  Then $X^m = I + TmX_1$ lies in the image of $j$.  Now suppose the result
holds for $r-1$ and consider the commutative diagram
$$\begin{array}{ccccccccc}
0  &  \ra  & K^{r-1}/K^r  &  \ra  & K/K^r &  \ra  & K/K^{r-1} & \ra & 1 \\
   &       & \downarrow   &    & \downarrow &     & \downarrow &  &     \\
0  &  \ra  &\U^{r-1}/\U^r &  \ra  &\U/\U^r & \ra  &\U/\U^{r-1} & \ra & 1. 
\end{array}$$  
Suppose $X \in \U/\U^r$.  Denote its image in $\U/\U^{r-1}$ by $\overline{X}$.  By the inductive hypothesis,
there is an integer $m\ne 0$ with $\overline{X}^m = \overline{Y}$ for some $\overline{Y} \in K/K^{r-1}$.  Choose
a lift $Y$ of $\overline{Y}$ in $K/K^r$.  Then $Y$ maps to $\overline{X}^m$ in $\U/\U^{r-1}$.  But $X^m$ also
maps to $\overline{X}^m$ so that $X^mY^{-1} = Z$ for some $Z\in \U^{r-1}/\U^r$.  Now, there exists some
$W \in K^{r-1}/K^r (\cong \lslnz)$ with $Z^p = W$ for some $p \ne 0$.  Since $Y=Z^{-1}X^m$, we have
\begin{eqnarray*}
Y^p & = & (Z^{-1}X^m)^p \\
    & = & Z^{-p}X^{mp} \qquad(\text{since}\; \U^{r-1}/\U^r\; \text{is central}) \\
    & = & W^{-1}X^{mp}.
\end{eqnarray*}
Thus, $X^{mp} = WY^p$ belongs to $K/K^r$ and the induction is complete.
\end{proof}

\begin{theorem}\label{malthm} The inclusion $K(R)\ra \U$ is the Malcev
completion.
\end{theorem}

\begin{proof} Since $\U= \varprojlim \U/\U^r$ and  since $\U/\U^r$ is the
Malcev completion of $K(R)/K^r(R)$, the theorem will follow
immediately if we can show that $K^r(R) = \Gamma^rK(R)$ for each $r$.
This follows from the next two lemmas.

\begin{lemma}[\cite{bousfield}, 13.6] \label{l1} Let $F^\bullet$ be a
central series in a group $G$ such that
\begin{enumerate}
\item The natural map $G\ra \varprojlim G/F^s$ is an isomorphism.
\item $F^s/F^{s+1}$ is torsion-free for $s\ge 1$.
\item The Lie product $G/F^2 \otimes F^s/F^{s+1} \ra F^{s+1}/F^{s+2}$
is surjective for $s\ge 1$.
\item $G/F^2$ is finitely generated.
\end{enumerate}
Then $\Gamma^sG = F^s$ for $s\ge 1$. \hfill $\qed$
\end{lemma}

\begin{lemma}[\cite{bousfield}, 13.4] \label{l2}  Let $G$ be a group
and denote by $\overline{G}$ the completion $\overline{G}=\varprojlim
G/\Gamma^rG$.  Then the following statements are equivalent.
\begin{enumerate}
\item The map $\overline{G} \ra \overline{\overline{G}}$ is an
isomorphism. 
\item The map $G/\Gamma^rG \ra \overline{G}/\Gamma^r\overline{G}$ is
an isomorphism for each $r\ge 1$. \hfill $\qed$
\end{enumerate}
\end{lemma}

\noindent {\em Completion of the proof of Theorem} \ref{malthm}.
Consider the group $$\overline{K} = \text{ker}(\powz
\stackrel{T=0}{\lra} \slnz)$$ with its $T$-adic filtration
$\overline{K}^\bullet$. Note that Lemma \ref{l1} shows that $\overline{K}^r =
\Gamma^r\overline{K}$ for each $r$:  the first two conditions are
clear, as is the fourth; the third condition follows since the Lie
algebra $\lslnz$ is perfect (it is here that we must exclude the case
$n=2$).  Since $$\overline{K} = \varprojlim K(R)/K^r(R)=
 \varprojlim K(R)/\Gamma^rK(R)$$
(the last equality is Corollary \ref{cor1}), and since
\begin{eqnarray*}
\overline{K} & = & \varprojlim \overline{K}/\overline{K}^r \\
             & = & \varprojlim \overline{K}/\Gamma^r\overline{K} \\
             & = & \overline{\overline{K}},
\end{eqnarray*} 
Lemma \ref{l2} implies that $K(R)/\Gamma^rK(R) \cong
\overline{K}/\Gamma^r\overline{K}$ for all $r$.  Consider the
commutative diagram
$$\begin{array}{ccc}
K(R)/\Gamma^rK(R) & \stackrel{\cong}{\ra} & \overline{K}/\Gamma^r\overline{K} \\
\downarrow       &                   &   \downarrow\cong \\
K(R)/K^r(R)   & \stackrel{\cong}{\ra}  &  \overline{K}/\overline{K}^r.
\end{array}$$ 
It follows that $K^r(R) = \Gamma^rK(R)$ and hence $\U$ is the Malcev
completion of $K(R)$.
\end{proof}

\begin{remark}  Even if Lemma \ref{l2} were not available, we could
still prove the result as follows.  Denote by $\M_r$ the Malcev
completion of $K(R)/\Gamma^rK(R)$, and by $\M = \varprojlim \M_r$ the
Malcev completion of $K(R)$.  Then the map $K(R)\ra \M$ factors
through $\overline{K}$.  Moreover, by the universal property of $\M$,
we get a unique map $\M\ra \U$ which is easily seen to be an
isomorphism since it has an inverse given by the universal property of
the Malcev completion $\overline{K} \ra \U$.
\end{remark}

\begin{cor}\label{h1} If $n\ge 3$, then $H_1(K(R),\zz) \cong
H_1(\overline{K},\zz) \cong \lslnz$.
\end{cor}

\begin{proof}  The first isomorphism follows from Lemma \ref{l2} and
the second isomorphism from Lemma \ref{l1}.
\end{proof}

\section{The Relative Completion of $\slnr$}\label{main}

We are now ready to prove the main result.

\begin{theorem}\label{mainthm}  If $n\ge 3$, then the map $\zt
\stackrel{t\mapsto T}{\lra} \powq$ (resp. $\zl \stackrel{t\mapsto
1+T}{\lra}\powq$) is the completion with respect to the map $\zt
\stackrel{t=0}{\lra} \slnq$ (resp. $\zl \stackrel{t=1}{\lra} \slnq$).
\end{theorem}

\begin{proof}  The relative completion is a proalgebraic group which
is an extension
\begin{equation}\label{ext}
1 \lra \cP \lra \G \lra \slnq \lra 1
\end{equation}
where $\cP$ is prounipotent.  By the universal property of $\U$, we
have a unique map $\Phi:\U \ra \cP$ induced by the map $K(R)\ra \cP$.
Since $H^1(\slnz,A)=0$ for all rational $\slnq$-modules $A$
\cite{rag}, we see that $\Phi$ is surjective (Proposition
\ref{surj}).  On the other hand, since $H_1(K(R),\zq) \cong \lslnq$ is finite
dimensional by Corollary \ref{h1}, since the action of $\slnz$ on
$H_1(K(R),\zq)$ extends to a rational representation of $\slnq$, and
since $H^2(\slnz,A)=0$ for all nontrivial rational representations $A$
of $\slnq$ \cite{borel}, Proposition \ref{inj} shows that
 the kernel of $\Phi$ is no larger than
$H_2(\slnz,\zq) = 0$.  Thus, $\U\cong \cP$,
and since the extension (\ref{ext}) is split (\cite[Prop. 4.4]{hain2}),
it follows that $\G \cong \powq$.
\end{proof}

\section{The Case $n=2$}\label{n2}

The proof of Theorem \ref{mainthm} breaks down in the case $n=2$ for a 
variety of reasons.

\begin{enumerate}
\item The Lie algebra ${\mathfrak sl}_2(\zz)$ is not perfect.
\item Raghunathan's theorem on the vanishing of $H^1(\slnz,A)$ does not
apply for $n=2$.
\item Borel's result for the vanishing of $H^2(\slnz,A)$ cannot be
strengthened to include $n=2$.
\end{enumerate}

However, one can make the following observations.  Denote by 
${\mathcal G}(\zz)$ the completion of $SL_2(\zz)$ relative to its
canonical inclusion in $SL_2(\zq)$, and by ${\mathcal G}(R)$ the completion
of $SL_2(R)$ ($R=\zz[t],\zz[t,t^{-1}]$) relative to the map
$SL_2(R)\stackrel{\mod {\mathfrak m}_R}{\lra} SL_2(\zq)$.  The group
${\mathcal G}(\zz)$ is {\em not} isomorphic to $SL_2(\zq)$; in fact, it
is an extension of $SL_2(\zq)$ by a free prounipotent group with
infinite dimensional $H_1$ (see \cite[Rmk. 3.9]{hain3}).

We have a commutative diagram
$$\begin{array}{ccccccccc}
1 & \lra & K(R) & \lra & SL_2(R) &\stackrel{\longleftarrow}{\lra} & SL_2(\zz) & \lra & 1 \\
  &      &\downarrow & & \downarrow &   & \downarrow &  &   \\
1 & \lra & {\mathcal W} & \lra & {\mathcal G}(R) &  
\stackrel{\longleftarrow}{\lra} & {\mathcal G}(\zz)
& \lra & 1.
\end{array}$$

\noindent The map $\Phi:{\mathcal G}(R)\ra {\mathcal G}(\zz)$ is induced
by the composition $SL_2(R)\ra SL_2(\zz)\ra {\mathcal G}(\zz)$ and
the map $\Psi:{\mathcal G}(\zz)\ra {\mathcal G}(R)$ is induced by the
composition $SL_2(\zz)\ra SL_2(R)\ra {\mathcal G}(R)$.  Since the 
composition $SL_2(\zz)\ra SL_2(R)\ra SL_2(\zz)$ is the identity, we
see that $\Phi \circ \Psi = \text{id}_{{\mathcal G}(\zz)}$.

If $n\ge 3$, then the map $\slnz\ra\slnq$ is the relative completion
so that the completion of $\slnr$ is an extension of the completion of
$\slnz$ by the Malcev completion of $K(R)$.  This leads us to make
the following conjecture.

\begin{conj}  The map $K(R)\lra {\mathcal W}$ is the Malcev completion.
\end{conj}

Note that there is at least some hope for this since ${\mathcal W}$ is
properly contained in the kernel of the map ${\mathcal G}(R)\ra SL_2(\zq)$
so that ${\mathcal W}$ is prounipotent.

\section{Cohomology}\label{coho}

In this section we provide evidence for the following conjecture.

\begin{conj}  If $n\ge 3$, then $H^2(\slnr,\zq) = 0$.
\end{conj}

Note that this conjecture is true for $n\ge 5$ for the following
reason.  If $n\ge 5$, then by van der Kallen's stability theorem
\cite{vdk}, we have
$$H_2(\slnr,\zz) \cong H_2(SL(R),\zz) \cong K_2(R).$$
It follows that if $n\ge 5$, then $H_2(\slnr,\zq) \cong K_2(R)\otimes
\zq$.  Since $K_2(\zz[t]) \cong K_2(\zz)$ and $K_2(\zz[t,t^{-1}])
\cong K_2(\zz) \oplus K_1(\zz)$, we see that $K_2(R)\otimes \zq = 0$.

The tool that we will use is continuous cohomology.  We define the
continuous cohomology of a group $\pi$ by
$$H^\bullet_{\text{cts}}(\pi,\zq) = \varinjlim
H^\bullet(\pi/\Gamma^r\pi,\zq).$$ 
The basic properties of continuous cohomology were established by Hain
\cite{hain1}.  We note the following facts.

\begin{prop}[\cite{hain1}, Thm. 5.1]  The natural map
$H^k_{\text{\em cts}}(\pi,\zq) \ra H^k(\pi,\zq)$ is an isomorphism for
$k=0,1$ and is injective for $k=2$. \hfill $\qed$
\end{prop}

The map on $H^2$ need not be surjective in general.  A group $\pi$ is
called {\em pseudo-nilpotent} if the natural map
$H^\bullet_{\text{cts}}(\pi,\zq) \ra H^\bullet(\pi,\zq)$ is an
isomorphism.  Examples of pseudo-nilpotent groups include the pure
braid groups, free groups and the fundamental groups of affine curves
over ${\mathbb C}$.

\begin{prop}[\cite{hain1}, Thm. 3.7] Let $\pi$ be a group with
$H_1(\pi,\zq)$ finite dimensional.  Let $\cP$ be the Malcev completion
of $\pi$ and denote by ${\mathfrak p}$ the Lie algebra of $\cP$.  Then
the natural map
$$H^\bullet_{\text{\em cts}}(\pi,\zq) \lra H^\bullet_{\text{\em cts}}({\mathfrak
p},\zq)$$
is an isomorphism.
\end{prop}

Thus, if $H_1(\pi,\zq)$ is finite dimensional, we can use Lie algebra
cohomology to obtain a lower bound on the dimension of
$H^2(\pi,\zq)$. We will not compute $H^2_{\text{cts}}(K(R),\zq)$ explicitly.  
However, we note the following result.

\begin{prop}  If $n\ge 3$, then $\dim H^2_{\text{\em cts}}(K(R),\zq) \ge 
(n^2-1)^2/4$.
\end{prop}

\begin{proof}  By a result of Lubotzky and Magid \cite{magid}, if $G$ is a 
nilpotent group with $b_1=\dim H_1(G,\zq)$ finite, then the second Betti
number $b_2$ satisfies $b_2\ge b_1^2/4$.  In the case of $K(R)/K^r(R)$, 
since $H_1(K/K^r,\zq) \cong \lslnq$, we see that $b_2(K/K^r) \ge 
(n^2-1)^2/4$ for each $r$.
\end{proof}

To show that $H^2(SL_n(R),\zq)$ vanishes, it would suffice to show the
following three things.

\begin{enumerate}
\item $H^2(\slnz,\zq)=0$.
\item $H^1(\slnz,H^1(K(R),\zq))=0$.
\item $H^0(\slnz,H^2(K(R),\zq))=0$.
\end{enumerate}

\noindent The first statement is clear.  The second follows from \cite{rag}
since $H^1(K(R),\zq)$ is the adjoint representation $\lslnq$.  The third 
statement is true for $n\ge 5$.

\begin{prop} If $n\ge 5$, then $H^0(\slnz,H^2(K(R),\zq))=0$.
\end{prop}

\begin{proof}  Consider the Hochschild--Serre spectral sequence
$$E_2^{p,q} = H^p(\slnz,H^q(K(R),\zq)) \Lra H^{p+q}(SL_n(R),\zq).$$
We know that $H^2(SL_n(R),\zq) = 0$ for $n\ge 5$ (see the remarks 
following Conjecture 5.1).  Note also that $H^2(\slnz,H^1(K(R),\zq)) = 0$
and $H^3(\slnz,\zq)=0$.  It follows that $d_2:E_2^{0,2}\ra E_2^{2,1}$
and $d_3:E_3^{0,2}\ra E_3^{3,0}$ are both the zero map and hence
$E_{\infty}^{0,2} = H^0(\slnz,H^2(K(R),\zq))$.  But this group must
vanish since $E_{\infty}^{1,1}$ and  $E_{\infty}^{2,0}$ do.
\end{proof}

The next result provides some evidence for the vanishing of the group
 $H^0(\slnz,H^2(K(R),\zq))$
when $n=3,4$.  We first state the following lemma, which can be proved via
direct computation.

\begin{lemma}\label{reps}  Let $\Gamma_{a_1,\dots ,a_{n-1}}$ be the 
irreducible $\slnq$-module with highest weight $(a_1+\cdots a_{n-1})L_1
+ \cdots + a_{n-1}L_{n-1}$, where $L_1,\dots , L_{n-1}$ are the
weights of the fundamental representation.  Then we have the following
isomorphisms of $\slnq$-modules:

\begin{enumerate}
\item ${\mathfrak sl}_3(\zq)\otimes {\mathfrak sl}_3(\zq) \cong \Gamma_{2,2} 
\oplus \Gamma_{3,0}\oplus \Gamma_{0,3}\oplus \Gamma_{1,1}\oplus \Gamma_{1,1}
\oplus \Gamma_{0,0}$
\item $\bigwedge^2 {\mathfrak sl}_3(\zq) \cong \Gamma_{3,0}\oplus \Gamma_{0,3}
\oplus \Gamma_{1,1}$
\item $\bigwedge^3 {\mathfrak sl}_3(\zq) \cong \Gamma_{2,2} \oplus \Gamma_{3,0}
\oplus \Gamma_{0,3} \oplus \Gamma_{1,1} \oplus \Gamma_{0,0}$
\item ${\mathfrak sl}_4(\zq)\otimes {\mathfrak sl}_4(\zq) \cong \Gamma_{2,0,2}
\oplus\Gamma_{2,1,0}\oplus\Gamma_{0,1,2}\oplus\Gamma_{0,2,0}\oplus
\Gamma_{1,0,1}\oplus\Gamma_{1,0,1}\oplus\Gamma_{0,0,0}$
\item $\bigwedge^2 {\mathfrak sl}_4(\zq) \cong \Gamma_{2,1,0}\oplus\Gamma_{0,1,2}
\oplus\Gamma_{1,0,1}$
\item 
$\bigwedge^3 {\mathfrak sl}_4(\zq)  \cong  \Gamma_{4,0,0} \oplus
\Gamma_{0,0,4} \oplus \Gamma_{1,2,1} \oplus \Gamma_{2,0,2} \oplus
\Gamma_{2,1,0} \oplus \Gamma_{0,1,2}$ \\
       $\mbox{} \qquad\qquad\qquad  \oplus \Gamma_{0,2,0} \oplus
\Gamma_{1,0,1} \oplus \Gamma_{0,0,0}$ \hfill $\qed$
\end{enumerate}
\end{lemma}

\begin{theorem}  If $n\ge 3$, then $H^0(\slnz,H^2_{\text{\em cts}}(K(R),\zq)) = 0$.
\end{theorem}

\begin{proof} We need only consider the cases $n=3,4$.  It suffices to
show that $H^0(\slnz,H^2(K/K^l,\zq))=0$ for each $l$.  We use Lie algebra
cohomology. Denote by ${\mathfrak u}$ the Lie algebra of $\U$ and consider
the $T$-adic filtration ${\mathfrak u}^\bullet$.  The Malcev Lie algebra of
$K/K^l$ is the Lie algebra ${\mathfrak u}_l = {\mathfrak u}/{\mathfrak u}^l$.
Observe that for each $l$, the quotient ${\mathfrak u}^{l-1}/{\mathfrak u}^l$
is isomorphic {\em as an $\slnq$-module} to the adjoint representation $\lslnq$,
but as a Lie algebra it is abelian ({\em i.e.,} to compute the bracket in
${\mathfrak u}^{l-1}/{\mathfrak u}^l$, we lift elements to ${\mathfrak u}^{l-1}$,
apply $[\quad ,\quad ]$, and project back; but the commutator of any two elements in
${\mathfrak u}^{l-1}$ lies in ${\mathfrak u}^l$ and so projects to 0).

We proceed by induction on $l$, beginning at $l=2$.  The Lie algebra 
${\mathfrak u}_2$ is abelian of dimension $n^2-1$; as an $\slnq$-module
it is the adjoint representation $\lslnq$.  Thus $H^2({\mathfrak u}_2,\zq)
\cong \bigwedge^2 \lslnq$ as an $\slnq$-module.  By Lemma \ref{reps}, parts
2 and 5, we see that $H^0(\slnz,H^2({\mathfrak u}_2,\zq))=0$.  Now, suppose
that $l>2$ and that $H^0(\slnz,H^2({\mathfrak u}_{l-1},\zq))=0$.
Consider the short exact sequence
$$0 \lra {\mathfrak u}^{l-1}/{\mathfrak u}^l \lra {\mathfrak u}_l
\lra {\mathfrak u}_{l-1} \lra 0.$$
The kernel is central in ${\mathfrak u}_l$.  Consider the Hochschild--Serre
spectral sequence
$$E_2^{p,q} = H^p({\mathfrak u}_{l-1},H^q({\mathfrak u}^{l-1}/
{\mathfrak u}^l,\zq)) \Lra H^{p+q}({\mathfrak u}_l,\zq).$$
We have isomorphisms of $\slnq$-modules:

\begin{enumerate}
\item $H^2({\mathfrak u}_{l-1},H^0({\mathfrak u}^{l-1}/{\mathfrak u}^l,\zq))
= H^2({\mathfrak u}_{l-1},\zq)$
\item $H^1({\mathfrak u}_{l-1},H^1({\mathfrak u}^{l-1}/{\mathfrak u}^l,\zq))
\cong \lslnq \otimes \lslnq$
\item $H^0({\mathfrak u}_{l-1},H^2({\mathfrak u}^{l-1}/{\mathfrak u}^l,\zq))
\cong \text{Hom}_\zq(\bigwedge^2\lslnq,\zq)$.
\end{enumerate}

\noindent By induction, the $\slnz$ invariants of the first module are
trivial and by Lemma \ref{reps}, parts 2 and 5, so are the invariants of
the last group.  It follows that $H^0(\slnz,E_{\infty}^{0,2}) = 0$.  Also,
since
$$H^1(\slnz,E_{\infty}^{0,1}) = H^1(\slnz,\lslnq) = 0,$$
the long exact cohomology sequence associated to the extension
$$0 \lra E_{\infty}^{0,1} \stackrel{d_2}{\lra} H^2({\mathfrak u}_{l-1},\zq)
\lra E_{\infty}^{2,0} \lra 0$$
shows that $H^0(\slnz,E_{\infty}^{2,0})=0$.  It remains to show that
 the group $H^0(\slnz,E_{\infty}^{1,1})$ vanishes.

Note that $E_2^{1,1}$ contains a copy of the trivial representation
(parts 1 and 4 of Lemma \ref{reps}).
However, the differential
$$d_2:E_2^{1,1} \lra H^3({\mathfrak u}_{l-1},\zq)$$
is easily seen to map the trivial representation onto a copy of the
trivial representation in the image (this copy arises from the
map in cohomology induced by the map ${\mathfrak u}_{l-1}\ra \lslnq$;
use parts 3 and 6 of Lemma \ref{reps}).
It follows that $E_{\infty}^{1,1}$ contains no copies of the trivial
representation and hence $H^0(\slnz,E_{\infty}^{1,1})=0$.  Thus
$$H^0(\slnz,H^2({\mathfrak u}_l,\zq)) = 0$$
and the induction is complete.
\end{proof}

One might conjecture that $K(R)$ is pseudo-nilpotent (we do not know
if this is the case).  If so, it would follow that $H^0(\slnz,H^2(K(R),\zq))=0$
and hence $H^2(SL_n(R),\zq)=0$ for $n\ge 3$.

\bibliographystyle{amsplain}

\end{document}